\documentclass[12pt,a4paper]{amsart}
\usepackage[T1]{fontenc}
\usepackage[utf8]{inputenc}
\usepackage{lmodern,microtype,amsmath,amssymb,amsthm,mathtools}
\usepackage[a4paper,left=23mm,right=23mm,top=23mm,bottom=24mm]{geometry}
\usepackage{booktabs,tabularx,array,enumitem,needspace}
\usepackage[hidelinks,breaklinks=true]{hyperref}
\setlist[itemize]{leftmargin=17pt,itemsep=3pt,topsep=3pt}
\setlist[enumerate]{leftmargin=20pt,itemsep=4pt,topsep=3pt}
\newcommand{\R}{\mathbb R}\newcommand{\C}{\mathbb C}
\newcommand{\dd}{\mathrm d}
\newcommand{\diag}{\operatorname{diag}}
\newcommand{\norm}[1]{\left\|#1\right\|}\newcommand{\Y}{\mathcal Y}

\newcolumntype{Y}{>{\raggedright\arraybackslash}X}
\newtheorem{theorem}{Theorem}[section]
\newtheorem{proposition}[theorem]{Proposition}
\newtheorem{lemma}[theorem]{Lemma}
\newtheorem{corollary}[theorem]{Corollary}
\theoremstyle{definition}
\theoremstyle{remark}\newtheorem{remark}[theorem]{Remark}
\hypersetup{pdftitle={Optimal boundary measurements and collision conditioning for biharmonic point defects},pdfauthor={Niyaz Tokmagambetov}}

\title[Optimal measurements for biharmonic defects]{Optimal measurement counts and collision conditioning for planar biharmonic point defects}
\author[N. Tokmagambetov]{Niyaz Tokmagambetov}
\address{Institute of Mathematics and Mathematical Modeling, 28 Shevchenko Street, 050010 Almaty, Kazakhstan}
\email{niyaz.tokmagambetov@gmail.com}
\date{}
\numberwithin{equation}{section}
\allowdisplaybreaks


\begin{document}

\begin{abstract}
We recover two planar biharmonic point defects from value and normal traces of the resolvent correction on an open circular arc.
For separated sites with inverse strengths of either sign bounded in magnitude above and away from zero, seven fixed real readings with a common input give globally Lipschitz recovery at a sufficiently small shift chosen from the admissible class.
Seven are necessary when the class interior contains a loop exchanging the sites.
For strengths in a compact positive interval, at any fixed positive shift, a finite measurement set independent of the minimum separation $h$ gives inverse losses $h^{-2}$ for positions and $h^{-3}$ for strengths.
If the admissible position set contains a disk, matching lower bounds hold even for the full nonlinear response.
Then noise $o(h^3)$ is necessary and sufficient for uniform position error $o(h)$ and strength error $o(1)$.
\end{abstract}
\maketitle
\par

\section{Introduction}

A localized perturbation of a plate can be described by its position and coupling strength. In the whole-plane model, exterior source and receiver fields can be chosen before the sites are known. For two defects, we seek a measurement count on separated configurations and a noise threshold for resolving individual sites as their separation tends to zero.

Let $H_0=\Delta^2$ on $L^2(\R^2)$ and consider the scalar form
\[
 q_{X,\alpha}[u]=\|\Delta u\|_2^2+
             \sum_{j=1}^N\alpha_j|u(x_j)|^2,
 \qquad u\in H^2(\R^2),\qquad \alpha_j\in\R\setminus\{0\}.
\]
The sites lie inside a fixed circle. At a resolvent point $-\lambda$,
with $\lambda>0$, we observe the defect contribution to the resolvent
kernel through its value and outward normal derivative on a fixed open
arc, in both source and observation variables. With
$b_j=\alpha_j^{-1}$, this response has the form
\[
 F_\lambda(X,b)=-W_XQ_\lambda(X,b)^{-1}W_X^*,\qquad
 Q_\lambda=\diag(b_j)+(G_\lambda(x_i-x_j))_{i,j}.
\]
Here $G_\lambda$ is the free Green function. Each column of $W_X$ contains its two observed traces for one site. Thus the measurement is an exterior response of the whole-space operator.

\subsection{Two main results and their parameter classes}

For two sites we prove the following results. All observations are fixed
from the a priori class, before the unknown configuration is chosen.

\begin{enumerate}
\item \emph{Seven global readings on a separated class.}
Suppose $x_1,x_2$ belong to a fixed compact interior set,
$|x_1-x_2|\ge d_0>0$, and
$0<b_*\le|b_j|\le b^*<\infty$. The inverse strengths may have
either sign. For every sufficiently small selected shift $\lambda>0$,
there are seven real source--test readings, with one common input, for
which
\[
 c_\lambda d_2(p,\widetilde p)
 \le |\mathcal M_\lambda(p)-\mathcal M_\lambda(\widetilde p)|
 \le C_\lambda d_2(p,\widetilde p).
\]
The distance $d_2$ includes both positions and inverse strengths,
modulo relabeling. Seven is necessary when the interior of the class contains an exchange
loop. The shift is chosen from the known class and then held fixed (Theorem~\ref{thm:seventwoglobal}).

\item \emph{Sharp collision losses for positive pairs.}
Fix any $\lambda>0$ and a positive interval
$0<b_-<b_+<\infty$. Let both configurations have separation at least
$h$ and inverse strengths in $[b_-,b_+]$. One finite boundary
compression $P_m$, independent of $h$, gives, for one matching $\pi$,
\[
 \|X-Y_\pi\|\le Ch^{-2}D_m(p,\widetilde p),\qquad
 \|b-\widetilde b_\pi\|\le Ch^{-3}D_m(p,\widetilde p),
 \quad D_m=\|P_m(F_\lambda(p)-F_\lambda(\widetilde p))P_m\|.
\]
If the position set contains a disk, both powers are optimal, even
for the complete response operator. The same powers hold for physical
strengths $\alpha_j=1/b_j$; see
Theorems~\ref{thm:sharptwoupper}--\ref{thm:sharptwolower}.
\end{enumerate}

Table~\ref{tab:classes} collects the assumptions used for the two results.
\begin{table}[htbp]
\centering
\caption{Measurement design and stability for two defects.}\label{tab:classes}
\begin{tabularx}{\linewidth}{@{}p{.21\linewidth}YY@{}}
\toprule
 & Seven-reading design & Collision estimates\\
\midrule
Inverse strengths & Either sign, magnitude in $[b_*,b^*]$ & Fixed positive interval\\
Separation & Fixed lower bound $d_0$ & Lower bound $h\downarrow0$\\
Spectral shift & Selected sufficiently small from the known class & Any fixed $\lambda>0$\\
Observations & Seven readings with one common input & One finite matrix independent of $h$\\
Optimality & Seven on classes with an interior exchange loop & Losses $h^{-2}$ and $h^{-3}$ when $K$ contains a disk\\
\bottomrule
\end{tabularx}
\end{table}
The constants include the norms of the prescribed boundary fields. In particular, the seven-reading count alone gives no bound on noise amplification as the selected shift tends to zero. Remark~\ref{rem:shiftsize} shows that for signed strengths the shift must be extremely small.

\subsection{Methods and previous work}

Vessella proved Lipschitz stability for recovering the locations and strengths of point sources from their generated potentials in three dimensions \cite{Vessella1992}. Recovery of operator point interactions from Green-kernel data has a
direct predecessor in Gesztesy and Ramm
\cite{Gesztesy2000Inhomogeneities}. The planar biharmonic point-extension
framework is classical \cite{Karpeshina1986Zero,Posilicano2008Extensions}.
For the present observations we first prove independence of Green
jets on every nonempty circle arc. Equality of response ranges then
identifies the sites, and a differential range argument identifies
infinitesimal changes of the positions and strengths.

The passage from injectivity and differential injectivity to
compact-class Lipschitz stability is developed by Bourgeois
\cite{Bourgeois2013Lipschitz}. Finite-measurement stability, including
manifold formulations and independent-data constructions, is treated in
\cite{Alberti2022Finite,Alberti2023Manifolds,Alberti2022Independent};
finite-electrode stability for a finite-dimensional conductivity class
is proved by Harrach \cite{Harrach2019Uniqueness}. We use this established
compactness mechanism in a supporting role. It does not determine the
optimal count or the separation powers proved here.

The seven-reading design combines two ingredients. First, seven real
moments determine an unordered measure with two planar nodes and
positive weights through the complex Prony formulas \cite{Batenkov2013Prony}.
Second, a small-shift argument makes the scaled physical response
charges positive even for signed inverse strengths. A Runge argument
realizes small perturbations of the reference moments by actual boundary
inputs and tests. The readings are exact response measurements; applying
the unperturbed Prony formulas to them is not asserted to be exact.
The lower count follows from orientation reversal when two
three-dimensional site--strength blocks are exchanged.

For Fourier spike data, Batenkov, Goldman and Yomdin
\cite{Batenkov2021Superresolution}, and Liu and Ammari for positive
sources \cite{Liu2023Positive}, establish different collision losses
for positions and amplitudes. Their two-source powers agree with the
powers $2$ and $3$ obtained here. Foundational sparse super-resolution
and generalized-spike conditioning are studied in
\cite{Donoho1992Superresolution,Batenkov2018Stability,Batenkov2023Generalized}.
For positive measures on the line, Denoyelle, Duval and Peyr\'e obtain support recovery by total-variation regularization at a noise and regularization scale set by $h^{2N-1}$ \cite{Denoyelle2017Support}. For two planar spikes at separation $h$, Poon and Peyr\'e show that if the noise and the regularization parameter are of order $h^4$ with bounded ratio, then total-variation regularization recovers exactly two spikes, with errors in the amplitudes and in the positions rescaled by $h$ bounded by a constant times the sum of these two quantities divided by $h^3$; they assume a nondegeneracy condition on the limiting precertificate \cite{Poon2019Sparse}. Both results concern additive observations. In the present model the internal Green matrix also varies with the defect positions.
Liu and Ammari also give a multidimensional Fourier resolution limit for distinguishing two sources from one \cite{Liu2024Resolution}. Here the number of sites is known and the recovery errors concern the individual positions and couplings.
The fixed-grid restriction in the minimax result of
\cite{Batenkov2023Generalized} and the additive Fourier observation
model are different from our off-grid nonlinear resolvent response.
Our upper estimate uses degree-three multivariate moments, controlled
uniformly over four arbitrarily clustered nodes by Kergin interpolation
\cite{Micchelli1980Constructive}. The lower estimate constructs physical
diagonal couplings whose response terms through degree two agree.
Thus the lower bound respects the internal Green matrix and concerns
the actual scalar strengths.

Convex recovery from Fourier data under separation conditions
\cite{Candes2014Theory} and sparse reconstruction for nonlinear
Foldy--Lax observations \cite{Alberti2024Localization} provide
complementary approaches. Point-source localization for biharmonic waves
is considered in \cite{Liu2026Identifying}. Neither additive source
recovery nor a linearization in the interaction strengths is substituted
for the full response used below.

Sections~\ref{sec:inv-separated}, \ref{sec:inv-seven} and~\ref{sec:inv-sharptwo} develop the trace model, the seven-reading construction and the collision bounds, respectively.

Parameter and measurement spaces are real. Fourier kernels and complex moments use their usual complexifications, with inner products linear in the second argument. We write $\|\cdot\|$ for the operator norm and $\asymp$ for two-sided bounds with constants independent of the separation.

\section{Boundary response and separated configurations}
\label{sec:inv-separated}

\subsection{Scalar defects and measurements}
\label{sec:inv-normalization}

Let $H_0=\Delta^2$ on $L^2(\R^2)$, with domain $H^4(\R^2)$.
Our Fourier convention is
\[
 \widehat f(\xi)=\int_{\R^2}e^{-ix\cdot\xi}f(x)\,\dd x,
 \qquad f(x)=(2\pi)^{-2}\int_{\R^2}
                 e^{ix\cdot\xi}\widehat f(\xi)\,\dd\xi.
\]
For $\lambda>0$, the Green function and its diagonal value are
\begin{equation}\label{eq:Bgreen}
 G_\lambda(x)=(2\pi)^{-2}\int_{\R^2}
             \frac{e^{ix\cdot\xi}}{|\xi|^4+\lambda}\,\dd\xi,
 \qquad g_\lambda=G_\lambda(0)=\frac1{8\sqrt\lambda}.
\end{equation}
For finitely many distinct sites and finite nonzero real strengths,
the form in the introduction is closed and semibounded on $H^2$.
Indeed, point values refer to the continuous Sobolev representative and
Fourier Cauchy--Schwarz gives
\[
 |u(x)|^2\le\frac1{8\sqrt t}
       \bigl(\|\Delta u\|_2^2+t\|u\|_2^2\bigr).
\]
At a sufficiently large $t$, the negative trace terms are a form
perturbation of relative bound less than one. The representation
theorem for closed semibounded forms defines the associated operator
$H_{X,b}$; see \cite{Kato1995Perturbation}.

Fix $R>0$, a nonempty open arc $\gamma$ of
$\Gamma=\partial B(0,R)$, and a nonempty compact
$K\Subset B(0,R)$ containing the sites. Use arc-length measure and the
outward normal of the disk, and set
\[
 \Y=L^2(\gamma;\R)\oplus L^2(\gamma;\R),\qquad
 \mathcal T u=(u,\partial_\nu u)|_\gamma.
\]
Initially $\mathcal T$ is defined on fields smooth near $\Gamma$.
Every Green field traced below has that regularity. Put
\[
 w(x)=\mathcal T G_\lambda(\cdot-x),\qquad
 h_a(x)=\langle a,w(x)\rangle_\Y \quad (a\in\Y).
\]
The column operator $W_X:\R^N\to\Y$ is given by $W_Xc=\sum_j c_jw(x_j)$.
At any $\lambda>0$ for which the real symmetric matrix
$Q_\lambda(X,b)$ is invertible,
\begin{equation}\label{eq:scalarphysical}
 \begin{gathered}
 Q_\lambda(X,b)=\diag(b_j)+M_\lambda(X),\qquad
 M_{\lambda,ij}(X)=G_\lambda(x_i-x_j),\\
 (H_{X,b}+\lambda)^{-1}=(H_0+\lambda)^{-1}
       -\Phi_XQ_\lambda^{-1}\Phi_X^*,\qquad
 \Phi_Xc=\sum_jc_jG_\lambda(\cdot-x_j),
 \end{gathered}
\end{equation}
and the bounded defect response is
\begin{equation}\label{eq:Bresponse}
             F_\lambda(X,b)=-W_XQ_\lambda(X,b)^{-1}W_X^*.
\end{equation}
To verify the resolvent identity, write
$u=R_0(\lambda)f-\Phi_Xc$, where $c_j=\alpha_ju(x_j)$, and
evaluate at the sites. This gives $Q_\lambda c=\Phi_X^*f$.
Invertibility of $Q_\lambda$ therefore also verifies that the declared
shift is a resolvent point. The correction kernel is smooth near
$\Gamma\times\Gamma$, so tracing both variables gives
\eqref{eq:Bresponse} without tracing an arbitrary $L^2$ operator.

An individual real reading is
$\langle a,F_\lambda(X,b)c\rangle_\Y$, with input $c$ and test
$a$. A fixed finite compression $P_mF_\lambda P_m$ records the entries
in a fixed orthonormal boundary basis. We use its operator norm; on a
fixed matrix size it is equivalent to the Euclidean norm of its
independent entries. A general bounded real linear functional of the
response need not be one individual source--test reading.
When the shift is fixed we also write $F=F_\lambda$.

\subsection{Qualitative and differential injectivity}
\label{sec:inv-channels}

\begin{lemma}[Independence on every fixed open arc]\label{lem:arc}
For any finite list of distinct site--multi-index pairs $(j,\alpha)$
with $|\alpha|\le3$ and distinct sites $x_j$ inside $\Gamma$, the fields
$\partial^\alpha G_\lambda(\cdot-x_j)$ have linearly independent two-traces on
$\gamma$. For the corresponding real column map $W$, one has
$W^*W>0$ and $\sigma(W)>0$, where
$\sigma(W)=\inf_{\norm c=1}\norm{Wc}_\Y$.
\end{lemma}
\begin{proof}
Partial fractions in $|\xi|^2$ and the modified Helmholtz kernel give
\[
 G_\lambda(r)=-\frac{\operatorname{Im} K_0(e^{i\pi/4}\lambda^{1/4}r)}
                          {2\pi\sqrt\lambda},\qquad r>0 .
\]
Indeed the heat integral gives $(2\pi)^{-1}K_0(mr)$ for real $m>0$.
Both the Fourier resolvent and this kernel continue holomorphically to $\operatorname{Re}m>0$ in tempered distributions. The power series and large-argument expansions for $K_0$ are given in \cite{NISTHandbook}.
The singular part starts with $r^2\log r/(8\pi)$. Thus third derivatives are locally integrable, and the boundary term from second derivatives on a circle of radius $r$ is $O(r(1+|\log r|))$. Derivatives through order three have no point-mass term. All these fields are analytic off the sites and decay exponentially with their exterior derivatives.

Let $u$ be a linear combination of the indicated jets with zero two-trace on the arc. Analyticity in the angle extends both zero traces to the entire circle. In the exterior $E=\R^2\setminus\overline{B(0,R)}$, twice integrating $(\Delta^2+\lambda)u=0$ against $u$ gives
\[
 \int_E (|\Delta u|^2+\lambda |u|^2)=0 .
\]
The clamped traces remove both boundary terms on $\Gamma$ and exponential decay removes the terms at infinity. Hence $u=0$ on $E$. Analytic continuation gives $u=0$ on the connected complement of the sites and local integrability gives $u=0$ as a distribution. Applying $\Delta^2+\lambda$ gives a zero combination of derivatives of point masses. Test functions with prescribed jets at each site show that all coefficients vanish. Compactness of the coefficient unit sphere gives the last assertion.
\end{proof}

We also use interior analytic regularity for distributional solutions of the homogeneous equation. Elliptic regularity and the analytic coefficients of $\Delta^2+\lambda$ give this directly \cite{Morrey1957Analytic}.

\begin{remark}[Why both traces are specified]
The proof uses the clamped exterior energy identity. It is not a proof
that one scalar trace on an arbitrary curve observes every jet. The circle
and its analyticity are used to pass from the observed arc to the whole
circle; the resulting stability constant depends on the chosen arc.
\end{remark}

\begin{proposition}[Scalar configurations and their linearizations]
\label{thm:sites}
At a fixed $\lambda>0$, the response of a finite scalar configuration
with distinct sites, finite nonzero strengths and invertible
$Q_\lambda$ determines its site count, sites and strengths up to
relabeling. The differential of this response map is injective on the
full $3N$-dimensional parameter space.
\end{proposition}
\begin{proof}
Set $B=Q_\lambda^{-1}$. Lemma~\ref{lem:arc} gives an injective
$W=W_X$, and $W^*$ is onto, so $\operatorname{ran}F=\operatorname{ran}W$.
For equal responses apply the lemma to the union of the two site lists.
Independence forces equality of their monopole spans only when the site
sets agree; their common rank is the number of sites. With these sites
fixed, put $L_W=(W^*W)^{-1}W^*$. Then
\begin{equation}\label{eq:Brecovery}
 B=-L_WFL_W^*,\qquad
 \diag(b_j)=B^{-1}-M_\lambda(X).
\end{equation}
This recovers the strengths. These inversions are used only at
admissible configurations; no arbitrary coefficient perturbation is
substituted for a physical scalar perturbation.

For the differential, let
$\dot p=(\xi_1,\ldots,\xi_N,t_1,\ldots,t_N)$ and
$\dot W_j=-\mathcal T(\xi_j\cdot\nabla G_\lambda(\cdot-x_j))$.
If $P$ is the orthogonal projection onto $\operatorname{ran}W$,
\[
 \dot F=-\dot WBW^*-WB\dot W^*+WB\dot QBW^*,\qquad
 (I-P)\dot FL_W^*=-(I-P)\dot WB.
\]
If $\dot F=0$, invertibility of $B$ makes every $\dot W_j$ belong to
$\operatorname{ran}W$. Independence of the monopoles and first jets
in Lemma~\ref{lem:arc} gives $\xi_j=0$ for every $j$. Thus
$\dot Q=\diag(t_j)$ and $\dot F=WB\diag(t_j)BW^*$. Left and right
inverses now give $t_j=0$ for every $j$.
\end{proof}

\subsection{Supporting finite-observation estimate}

Fix $N$, a compact $K\Subset B(0,R)$, $d_0>0$, $0<b_*\le b^*<\infty$
and $\kappa>0$. Let $\mathcal P_N$ be any nonempty closed, permutation-invariant
class of scalar configurations satisfying
\begin{equation}\label{eq:separatedclass}
 x_j\in K,\quad |x_i-x_j|\ge d_0\ (i\ne j),\quad
 b_*\le |b_j|\le b^*,\quad
 \sigma_{\min}(Q_\lambda(X,b))\ge\kappa,
\end{equation}
where $Q_\lambda=\diag(b_j)+M_\lambda(X)$ and
$M_{\lambda,ij}=G_\lambda(x_i-x_j)$. Parameters have fixed units. Use
the relabeling distance
\begin{equation}\label{eq:permdistance}
 d_N(p,q)=\min_{\pi\in\mathfrak S_N}
 \left(\sum_{j=1}^N\bigl(|x_j-y_{\pi(j)}|^2+
                   |b_j-\widetilde b_{\pi(j)}|^2\bigr)\right)^{1/2}.
\end{equation}
The class is compact. The lower bound on $|b_j|$ excludes undefined
finite strengths, and the last condition is a spectral gap at the
declared measurement shift. Neither condition supplies a separation-free
stability constant.

\begin{proposition}[Separated compact classes]
\label{thm:unknownfinite}
On $\mathcal P_N$, the response $F=F_\lambda$ satisfies
\begin{equation}\label{eq:unknownLip}
 c_Nd_N(p,q)\le\|F(p)-F(q)\|\le C_Nd_N(p,q).
\end{equation}
For any complete real orthonormal boundary basis, a sufficiently large
fixed projection $P_m$ preserves this estimate, with lower constant
$c_N/2$. The constants and $m$ depend on the known separated class.
The statement also holds in the physical strengths $\alpha_j=1/b_j$.
\end{proposition}
\begin{proof}
On the open set of distinct sites, nonzero $b_j$ and invertible $Q_\lambda$, the response is smooth. The permutation action is free, so a local labeling gives a smooth chart of the quotient in which $d_N$ is Euclidean near the diagonal. Proposition~\ref{thm:sites} gives injectivity and an injective differential. The compact-class Lipschitz theorem \cite{Bourgeois2013Lipschitz,Alberti2023Manifolds} gives \eqref{eq:unknownLip} in these charts, and compactness and injectivity separate the remaining pairs.

For the finite observations apply the compression principle of \cite{Alberti2022Finite,Alberti2023Manifolds} to $T\mapsto P_mTP_m$. To see its hypotheses here, the normalized secants $(F(p)-F(q))/d_N(p,q)$ have compact closure. Along a convergent sequence of pairs with coincident limits, local labeling and Taylor expansion reduce this to compactness of the unit sphere in $\R^{3N}$. Distinct limits follow by continuity. All limits are compact operators, on which the contractive compressions converge in norm, uniformly on this compact closure. An error of at most $c_N/2$ gives the stated bound. Finally $|1/b-1/\widetilde b|=|b-\widetilde b|/|b\widetilde b|$ transfers the estimates to physical strengths of either sign.
\end{proof}

\section{An optimal seven-reading design for separated pairs}
\label{sec:inv-seven}\label{sec:seventwo}

We first separate local dimension from global unordered geometry. The
six coordinates of a labeled pair do not by themselves give a global
six-reading theorem for unordered pairs. After the lower bound we
construct seven exact source--test readings at one designable shift.

\subsection{Local dimension and the exchange obstruction}

\begin{theorem}[Exactly $3N$ local source--test readings]
\label{thm:localoptimal}
At any admissible $N$-site parameter $p_0$, there are $3N$ real
boundary vectors $a_1,\ldots,a_{3N}\in\Y$ and a neighborhood
$\mathcal U$ of $p_0$, modulo relabeling, such that
\[
 p\longmapsto\big(\langle a_j,F(p)a_j\rangle\big)_{j=1}^{3N}
\]
is bi-Lipschitz on a sufficiently small compact parameter ball in
$\mathcal U$. The vectors are fixed from the nominal parameter $p_0$.
No smaller number of scalar readings with a Lipschitz inverse works
on an open parameter patch. Thus $3N$ is the optimal local scalar
count, and these readings are individual source--test experiments.
\end{theorem}
\begin{proof}
Proposition~\ref{thm:sites} proves that $DF(p_0)$ is injective.
Its range is a $3N$-dimensional real vector space of self-adjoint
operators. Quadratic functionals $T\mapsto\langle a,Ta\rangle$
separate this space: if all vanish, real polarization gives $T=0$.
Their restrictions therefore span its dual. Choose $3N$ independent
restrictions. The derivative of the displayed measured map is an
invertible square matrix in any local labeling chart. The inverse
function theorem applies. More explicitly, on a small convex ball
the derivative differs from its central matrix by less than half
that matrix's least singular value; integration on line segments
gives the two Lipschitz bounds. Local chart distance is equivalent
to $d_N$. For necessity, differentiate a lower Lipschitz estimate at an
interior point of a local labeling chart. An injective differential
from $\mathbb R^{3N}$ to $\mathbb R^m$ requires $m\ge3N$.
\end{proof}

\begin{proposition}[An exchange-loop obstruction for unordered sites]
\label{prop:topologicalcount}
For every $N\ge2$ there are compact separated positive-strength
classes of the form \eqref{eq:separatedclass} for which any
continuous injective map into $\R^m$ requires $m\ge3N+1$. More generally, the same bound holds whenever the
interior of the admissible class in the full $3N$-dimensional quotient manifold contains a loop whose lift exchanges two site--strength blocks and fixes the others.
\end{proposition}
\begin{proof}
The ordered interior parameter manifold has its usual orientation
in $3N$ real coordinates, grouped as $(x_j^{(1)},x_j^{(2)},b_j)$.
Permuting two blocks changes orientation by $(-1)^{3\cdot3}=-1$.
Since the permutation action is free, an exchange loop in the
quotient reverses its local orientation. A neighborhood of that
loop is therefore nonorientable.

Put $n=3N$ and let $U$ be an open neighborhood of the loop in the interior of the quotient. If $m\le n$, compose a continuous injective map $U\to\R^m$ with the coordinate inclusion into $\R^n$. Invariance of domain makes the composite a homeomorphism onto an open subset of $\R^n$ \cite{Bredon1993Topology}. For $m<n$ its image lies in a subspace with empty interior. For $m=n$ the open image is orientable, whereas transport of a local orientation around the loop reverses its sign. Orientability is invariant under homeomorphisms, so both possibilities are excluded.

To obtain a concrete compact class, place two sites at opposite
points of a small interior circle, give them equal positive inverse
strengths, and rotate the pair through a half-turn. Keep all other
sites fixed, separated from that circle and each other. The endpoint
is the starting unordered configuration, and the ordered lift is a
transposition. Choose a small compact neighborhood of this loop,
with strict separation, positions away from the observation circle,
and $b_j$ in a positive interval. This neighborhood contains an open
neighborhood of the loop. The free Green matrix is positive
semidefinite by its Fourier Gram representation, so
$Q_\lambda\ge b_*I$; a uniform spectral gap follows. This supplies
the asserted admissible class. The obstruction concerns global
unordered classes; it does not affect Theorem~\ref{thm:localoptimal}.
\end{proof}

\subsection{Boundary approximation and seven reference moments}

\begin{lemma}[Interior approximation from the observation arc]
\label{lem:arcrunge}
Let $D_1\Subset D_2\Subset B(0,R)$ be concentric open disks. Every
real solution of $(\Delta^2+\lambda)h=0$ on a neighborhood of
$\overline D_2$ is a $C^1(\overline D_1)$ limit of functions $h_a$.
\end{lemma}
\begin{proof}
Suppose $f\in L^2(D_2)$ annihilates every $h_a$, and extend $f$ by
zero to the plane. Then $u=R_0(\lambda)f\in H^4(\R^2)$ satisfies
$\mathcal Tu=0$ on $\gamma$, by Fubini and the adjoint identity.
The source is separated from the observation circle. Both circle
traces are consequently real analytic and vanish on the whole circle.
The exterior integration argument of Lemma~\ref{lem:arc} gives $u=0$
outside that circle. Analytic continuation for the homogeneous
constant-coefficient equation on the connected set
$\R^2\setminus\overline D_2$ then gives
$\operatorname{supp}u\subset\overline D_2$.

Choose a smooth compactly supported cutoff $\chi$ in the neighborhood
where $h$ is defined, equal to one near $\overline D_2$. Distributional
integration by parts yields
\[
 \int_{D_2}fh=\langle(\Delta^2+\lambda)u,\chi h\rangle
             =\langle u,(\Delta^2+\lambda)(\chi h)\rangle=0.
\]
The last factor vanishes near the support of $u$. The orthogonal
complement criterion in the real Hilbert space $L^2(D_2)$ proves
that $h$ lies in the closure of the restrictions of the $h_a$.
All approximation errors solve the homogeneous equation in $D_2$.
The interior elliptic estimate, for example its $H^4(D_1')$ version
with $D_1\Subset D_1'\Subset D_2$, and Sobolev embedding (see \cite{Gazzola2010Polyharmonic}) upgrade
$L^2(D_2)$ convergence to $C^1(\overline D_1)$ convergence.
\end{proof}

Identify $x=(x^{(1)},x^{(2)})$ with $z=x^{(1)}+ix^{(2)}$. The
seven real functions used as reference tests are
\begin{equation}\label{eq:sevenpolynomials}
 (P_0,\ldots,P_6)=(1,\operatorname{Re}z,\operatorname{Im}z,
 \operatorname{Re}z^2,\operatorname{Im}z^2,
 \operatorname{Re}z^3,\operatorname{Im}z^3).
\end{equation}
Their inversion is the classical two-node complex Prony calculation;
see also \cite{Batenkov2013Prony}. We include its stability
argument and separately realize a perturbation of these tests through
the physical two-trace response.

\begin{lemma}[Seven moments and stability under perturbing the tests]
\label{lem:sevenmoments}
Fix $r,d_0>0$ and $0<c_-<c_+<\infty$. On the class of unordered
two-atom measures
\[
 \mu=A_1\delta_{z_1}+A_2\delta_{z_2},\qquad
 |z_i|\le r,\quad |z_1-z_2|\ge d_0,\quad A_i\in[c_-,c_+],
\]
the seven real moments $\int P_j\,d\mu$, $0\le j\le6$, form
a bi-Lipschitz map for the Euclidean relabeling metric on nodes and
weights. This remains true if the $P_j$ are replaced by any real
$C^1$ functions $h_j$ sufficiently close to them in
$C^1(\overline{B(0,r)})$, with a uniform closeness threshold.
\end{lemma}
\begin{proof}
Combine the real readings into $m_0\in\R$ and $m_1,m_2,m_3\in\C$,
where $m_k=A_1z_1^k+A_2z_2^k$. Then
\begin{equation}\label{eq:twopronyformulas}
 \begin{gathered}
 \Delta_m=m_0m_2-m_1^2=A_1A_2(z_1-z_2)^2,\\
 S=\frac{m_0m_3-m_1m_2}{\Delta_m}=z_1+z_2,
 \qquad P=\frac{m_1m_3-m_2^2}{\Delta_m}=z_1z_2.
 \end{gathered}
\end{equation}
The nodes are the unordered roots of $Z^2-SZ+P$, and
\[
 A_1=\frac{m_1-m_0z_2}{z_1-z_2},\qquad A_2=m_0-A_1.
\]
In particular the map is injective. The denominator satisfies
$|\Delta_m|\ge c_-^2d_0^2$. On the compact parameter class the
coefficients $S,P$ are therefore uniformly Lipschitz functions of
the moments. For two nearby moment vectors, evaluating the second
root polynomial at a root of the first gives
\[
 |(z_i-\widetilde z_1)(z_i-\widetilde z_2)|
 \le C|m-\widetilde m|.
\]
For sufficiently small moment difference each root is within $d_0/3$
of one new root. Separation forces a bijective matching, and the other
factor is at least $2d_0/3$. The matched root errors are consequently
at most $C'|m-\widetilde m|$. The weight formula has the same
Lipschitz property. For larger moment differences, boundedness of the
parameter class gives the lower Lipschitz estimate after decreasing its
constant. The upper estimate follows directly from the moment formula.

Finally put $e_j=h_j-P_j$. For a minimizing matching of the two measures,
bounded weights and integration on segments in the disk give
\[
 \left|\sum_i A_i e_j(z_i)-
              \sum_i\widetilde A_i e_j(\widetilde z_i)\right|
 \le C\|e_j\|_{C^1}\,
       d\big((z,A),(\widetilde z,\widetilde A)\big).
\]
If all seven errors are sufficiently small,
the Lipschitz norm of the perturbation is less than half the lower
constant just proved. Subtracting it proves both bounds for the
perturbed test map. Exact Prony formulas are asserted for the reference
moments, not for the perturbed physical readings.
\end{proof}

\begin{lemma}[Small-shift fields approximating the seven tests]
\label{lem:smallshiftfields}
Fix $r<R$. For every sufficiently small $\lambda>0$, there are real
boundary vectors $a_0,\ldots,a_6\in\Y$ whose fields
$h_j(x)=\langle a_j,\mathcal T G_\lambda(\cdot-x)\rangle_\Y$
satisfy
\begin{equation}\label{eq:smallshiftfields}
 \max_{0\le j\le6}\|h_j-P_j\|_{C^1(\overline{B(0,r)})}
                                      \le C_r\lambda.
\end{equation}
In particular $h_0=1+O(\lambda)$ uniformly. The vectors may depend
on the selected shift and the known disk, but not on the unknown sites.
\end{lemma}
\begin{proof}
If $P$ is a harmonic homogeneous polynomial of degree $k$, then
\[
 \Delta(|x|^{2s}P)=4s(s+k)|x|^{2s-2}P.
\]
Define $a_{k,0}=1$ and
\[
 a_{k,j+1}=\frac{-\lambda a_{k,j}}
 {16(2j+2)(2j+k+2)(2j+1)(2j+k+1)},\qquad
 H_{P,\lambda}(x)=P(x)\sum_{j=0}^{\infty}a_{k,j}|x|^{4j}.
\]
The factorial growth of the denominator gives convergence with every
derivative on every bounded disk. Applying the displayed Laplacian
identity twice proves $(\Delta^2+\lambda)H_{P,\lambda}=0$.
On each fixed disk and for $0<\lambda\le1$,
$\|H_{P,\lambda}-P\|_{C^1}\le C\lambda$.
Apply Lemma~\ref{lem:arcrunge} between two slightly larger disks
inside $B(0,R)$ to the seven exact solutions $H_{P_j,\lambda}$
associated with the real polynomials in \eqref{eq:sevenpolynomials},
choosing the $C^1$ approximation error at most $\lambda$.
These exact boundary-generated fields prove
\eqref{eq:smallshiftfields}. No bound uniform in $\lambda$ is
claimed for the boundary-vector norms.
\end{proof}

\begin{lemma}[Uniformly positive scaled charges at a small shift]
\label{lem:scaledcharges}
Fix $K\Subset B(0,r)$, $d_0>0$ and $B<\infty$. Suppose
$x_1,x_2\in K$, $|x_1-x_2|\ge d_0$, and $|b_i|\le B$.
With $h_0$ from Lemma~\ref{lem:smallshiftfields}, put
\[
 Q_\lambda=\begin{pmatrix}b_1+g_\lambda&G_\lambda(x_1-x_2)\\
 G_\lambda(x_1-x_2)&b_2+g_\lambda\end{pmatrix},\qquad
 q=Q_\lambda^{-1}(h_0(x_1),h_0(x_2))^T.
\]
For all sufficiently small $\lambda>0$, uniformly on this class,
$Q_\lambda$ is positive definite and
\begin{equation}\label{eq:scaledcharges}
                       g_\lambda q_i\longrightarrow\tfrac12
                       \qquad(i=1,2).
\end{equation}
In particular $A_i=g_\lambda q_i\in[1/4,3/4]$ after decreasing
the shift threshold. This holds for either sign of the physical
inverse strengths.
\end{lemma}
\begin{proof}
Write $\rho=|x_1-x_2|$ and
$D_\lambda(\rho)=g_\lambda-G_\lambda(x_1-x_2)$. Fourier inversion gives
\[
 D_\lambda(\rho)=(2\pi)^{-2}\int_{\R^2}
       \frac{1-\cos(\xi\cdot(x_1-x_2))}{|\xi|^4+\lambda}\,d\xi.
\]
The numerator is nonnegative. For $|\xi|\le\rho_0$, where
$\rho_0(2r)\le1$, use $1-\cos t\ge t^2/4$; its angular integral
is a positive constant times $\rho^2|\xi|^2$. For the upper bound
use $1-\cos t\le t^2/2$ up to a fixed frequency and $1-\cos t\le2$
above it. Integrating $s^3/(s^4+\lambda)$ on the low-frequency
interval proves the uniform estimates
\begin{equation}\label{eq:greenlowshiftgap}
 c d_0^2\log(1/\lambda)-C
 \le D_\lambda(\rho)\le C_r(1+\log(1/\lambda))
                    \quad(d_0\le\rho\le2r).
\end{equation}
Thus $D_\lambda\to\infty$ and $D_\lambda/g_\lambda\to0$
uniformly, since $g_\lambda=(8\sqrt\lambda)^{-1}$.
Put $G=g_\lambda-D_\lambda>0$ and $T_i=b_i+D_\lambda>0$.
Then
\[
 Q_\lambda=\diag(T_1,T_2)+G\boldsymbol1\boldsymbol1^T>0,
 \qquad \det Q_\lambda=T_1T_2+G(T_1+T_2).
\]
The first charge is
\[
 g_\lambda q_1=
 \frac{g_\lambda[T_2h_0(x_1)+G(h_0(x_1)-h_0(x_2))]}
      {T_1T_2+G(T_1+T_2)}.
\]
Here $T_i/D_\lambda\to1$, $G/g_\lambda\to1$,
$T_1T_2/(g_\lambda D_\lambda)\to0$, and $h_0=1+O(\lambda)$.
The second numerator term, divided by $D_\lambda$, is bounded by
$C g_\lambda\lambda/D_\lambda\to0$.
Dividing numerator and denominator by $g_\lambda D_\lambda$
therefore gives $1/2$. The second charge is identical with indices
interchanged. Every limit used here is uniform under the declared
a priori bounds. Positivity of $Q_\lambda$ is a statement about the
finite response matrix; it does not assert positivity of the perturbed
plate operator for signed couplings.
\end{proof}

\begin{theorem}[Seven global scalar readings, with a designable shift]
\label{thm:seventwoglobal}
Fix $K\Subset B(0,R)$, $d_0>0$, and
$0<b_*\le b^*<\infty$. Let $\mathcal C_2$ be any nonempty closed unordered
class satisfying
\[
 x_i\in K,\qquad |x_1-x_2|\ge d_0,\qquad
                       b_*\le|b_i|\le b^*.
\]
There is $\lambda_0>0$, depending only on these bounds, such that
for every selected $0<\lambda<\lambda_0$, seven fixed real
source--test readings at the one energy $-\lambda$ give a global
Lipschitz inverse for both positions and strengths. They have one
common input $a_0$ and tests $a_0,\ldots,a_6$ on the prescribed open
arc. More precisely,
\begin{equation}\label{eq:sevenphysicaldata}
 \mathcal M_\lambda(p)=
 \bigl(-g_\lambda\langle a_j,F_\lambda(p)a_0\rangle_\Y\bigr)_{j=0}^6,
 \qquad
 c_\lambda d_2(p,\widetilde p)\le
 |\mathcal M_\lambda(p)-\mathcal M_\lambda(\widetilde p)|
 \le C_\lambda d_2(p,\widetilde p).
\end{equation}
If the interior of the class contains an exchange loop as in
Proposition~\ref{prop:topologicalcount}, seven is necessary and
sufficient, and the lower bound holds for every continuous injective
measurement map on that class. Consequently seven is the optimal universal scalar
count for compact separated two-site classes when the measurement
shift can be chosen from their a priori bounds.
\end{theorem}
\begin{proof}
Choose $r<R$ with $K\Subset B(0,r)$. By
Lemmas~\ref{lem:smallshiftfields} and~\ref{lem:scaledcharges}, at every
sufficiently small selected shift the fields $h_j$ are $C^1$ close to
$P_j$ and the scaled charges $A_i=g_\lambda q_i$ lie in
$[1/4,3/4]$. The physical response formula gives the exact identity
\begin{equation}\label{eq:sevenchargedata}
 -g_\lambda\langle a_j,F_\lambda(X,b)a_0\rangle
                         =\sum_{i=1}^2 A_i h_j(x_i).
\end{equation}
The smallness threshold in Lemma~\ref{lem:sevenmoments}, applied on
the full class with these fixed weight bounds and separation $d_0$,
is independent of $\lambda$. Decrease $\lambda_0$ so that
\eqref{eq:smallshiftfields} lies below it. The seven exact physical
readings then recover the unordered nodes and scaled charges with
a global Lipschitz inverse.

For fixed $\lambda$, passage between $(X,b)$ and $(X,A)$ is itself
bi-Lipschitz on the physical class. In the forward direction use
$q=Q_\lambda^{-1}h_0(X)$, the uniform bound on $Q_\lambda^{-1}$,
and $Q^{-1}-\widetilde Q^{-1}=Q^{-1}(\widetilde Q-Q)\widetilde Q^{-1}$.
In the reverse direction the exact force balance gives
\begin{equation}\label{eq:sevenstrengthrecovery}
 b_i=\frac{g_\lambda h_0(x_i)
                 -G_\lambda(x_i-x_{3-i})A_{3-i}}{A_i}-g_\lambda,
                         \qquad i=1,2.
\end{equation}
The denominators are at least $1/4$, and all functions and first
derivatives are bounded on the fixed compact position class, with
positive separation. These formulas give uniform Lipschitz estimates
for each matching; minimization over the two matchings yields the
quotient estimates. The same applies to $\alpha_i=1/b_i$ because
$|b_i|\ge b_*$. Combining the estimates proves
\eqref{eq:sevenphysicaldata}. The exact reconstruction may be defined
as inversion of the perturbed seven-test map followed by
\eqref{eq:sevenstrengthrecovery}; applying the reference Prony
formulas directly to uncorrected physical data is not asserted exact.

An exchange loop in the interior requires at least $3\cdot2+1=7$ scalar readings
by Proposition~\ref{prop:topologicalcount}. Such compact classes exist
with positive strengths, so six cannot be a universal count.
The construction supplies seven for every declared two-site class,
including signed ones. The lower bound need not apply to a restricted
class that admits a global labeling, and no claim is made that every
such restricted class needs seven.
\end{proof}

\begin{corollary}[Noisy seven-reading reconstruction]
\label{cor:sevennoise}
At the selected shift, suppose
$|y-\mathcal M_\lambda(p)|\le\eta$ in $\mathbb R^7$.
A residual minimizer $\widehat p$ over $\mathcal C_2$ exists and
satisfies $d_2(\widehat p,p)\le2\eta/c_\lambda$.
\end{corollary}
\begin{proof}
The class is compact and the measured map is continuous. Its minimizer
has residual no larger than the true parameter, so the two exact
measurement vectors differ by at most $2\eta$. Apply
\eqref{eq:sevenphysicaldata}.
\end{proof}

If the response error has operator norm at most $\eta_{\rm op}$, the error in the vector of seven readings is at most
\[
 g_\lambda\|a_0\|\left(\sum_{j=0}^6\|a_j\|^2\right)^{1/2}\eta_{\rm op}.
\]
This factor records the cost of the boundary fields at the selected shift.

\begin{remark}[Size of the selected shift]\label{rem:shiftsize}
For signed strengths the threshold $\lambda_0$ in Theorem~\ref{thm:seventwoglobal}
is extremely small. In the notation of the proof of Lemma~\ref{lem:scaledcharges},
$g_\lambda q_1=T_2/(T_1+T_2)+o(1)$ and $g_\lambda q_2=T_1/(T_1+T_2)+o(1)$ as
$\lambda\to0$. For $b_1=-b^*$ and $b_2=b^*$ these equal
${\tfrac12\bigl(1\pm b^*/D_\lambda(\rho)\bigr)}+o(1)$. Positive charges therefore
require $D_\lambda(\rho)>b^*$, and the window $[1/4,3/4]$ requires
$D_\lambda(\rho)\ge2b^*$. Evaluating the low-frequency integral behind
\eqref{eq:greenlowshiftgap} exactly gives
\[
 D_\lambda(\rho)=\frac{\rho^2}{32\pi}\log\frac1\lambda+O(1),
 \qquad d_0\le\rho\le2r,
\]
so $\log(1/\lambda_0)$ is of order $b^*/d_0^2$. For $b^*=1$ and $d_0=1/2$,
positive charges need $\lambda<10^{-171}$ and the window needs
$\lambda<10^{-346}$; at the latter shift $g_\lambda\approx10^{172}$. When both
strengths are positive, $T_i\ge b_*$ and $g_\lambda q_i$ lies, up to $o(1)$, in
$[b_*/(b_*+b^*),\,b^*/(b_*+b^*)]$ for every small shift. Using this window in
place of $[1/4,3/4]$ avoids the logarithmic condition.
\end{remark}

\section{Sharp collision conditioning for positive pairs}\label{sec:inv-sharptwo}

\label{sec:sharptwo}
We now let the separation vary while keeping the physical strengths in
a fixed positive interval. The estimates below concern the full nonlinear
response, including its internal Green matrix. They distinguish the
position error from the strength error. Fix
\[
 0<b_-<b_+<\infty,\qquad
 \mathcal P_2(h)=\{(X,b):x_1,x_2\in K,\ |x_1-x_2|\ge h,
                         \ b_j\in[b_-,b_+]\}/\mathfrak S_2.
\]
Write $w(x)=\mathcal T G_\lambda(\cdot-x)$,
$Q=\diag(b_1,b_2)+(G_\lambda(x_i-x_j))$, and
$F(X,b)=-W_XQ^{-1}W_X^*$. The Fourier Gram formula gives
$Q\ge b_-I$, also at coincident sites. All constants in this section
may depend on $K,\lambda,\gamma,b_-,b_+$ and the fixed trace norm,
but not on $h$. Vector and tensor norms below are Euclidean; their
finite-dimensional equivalence only changes these constants.

\subsection{Uniform reduction to four-node moments}
\begin{lemma}[Four-node moment comparison]\label{lem:fourmoment}
There are $\rho>0$ and $c_0>0$, uniform for centres $c$ in a compact
interior set, such that, for at most four nodes $z_i\in B(c,\rho)$
and real coefficients $a_i$,
\begin{equation}\label{eq:fourmoment}
 \left\|\sum_i a_iw(z_i)\right\|_\Y
 \ge c_0\left(\sum_{|\alpha|\le3}
       \left|\sum_i a_i(z_i-c)^\alpha\right|^2\right)^{1/2}.
\end{equation}
For any fixed complete real orthonormal boundary basis and its nested
coordinate projections $P_m$, the same constants, after decreasing
$c_0$ if necessary, work with $P_mw$ in place of $w$ for every
sufficiently large $m$. There is no lower bound on the mutual
distances of the four nodes.
\end{lemma}
\begin{proof}
The map $w$ is smooth in an interior neighbourhood, with values in
$\Y$. By Lemma~\ref{lem:arc}, the ten-column jet map
\[
 J_c v=\sum_{|\alpha|\le3}
           \frac{\partial_x^\alpha w(c)}{\alpha!}v_\alpha
\]
is injective. Its least singular value has a positive lower bound
as $c$ ranges over a compact interior set.

We use the degree-three simplex interpolation formula, a special case
of Kergin interpolation; see Micchelli~\cite{Micchelli1980Constructive}. To make its
uniformity at arbitrary, possibly repeated nodes explicit, put
\[
 \int_{[v_0,\ldots,v_k]} f
 =\int_{\substack{t_j\ge0\\t_1+\cdots+t_k\le1}}
 f\left(v_0+\sum_{j=1}^k t_j(v_j-v_0)\right)\dd t_1\cdots\dd t_k.
\]
These are parameter-simplex integrals, of mass $1/k!$, rather than
geometric volume integrals. Pad a shorter node list by repetitions
and zero coefficients and write its four nodes as $z_0,\ldots,z_3$.
The $\Y$-valued polynomial
\begin{equation}\label{eq:simplexpoly}
 I_Zw(x)=w(z_0)+\sum_{k=1}^3
 \int_{[z_0,\ldots,z_k]}
 D^kw(\xi)[x-z_0,\ldots,x-z_{k-1}]
\end{equation}
interpolates $w$ at every $z_i$. These are Bochner integrals. Applying each continuous linear functional on $\Y$ to the formula recovers scalar Kergin interpolation \cite{Micchelli1980Constructive}, including by continuity at repeated nodes. The functionals separate vectors, so the interpolation identity also holds in $\Y$.

Expand \eqref{eq:simplexpoly} in the monomials $(x-c)^\alpha$.
Let $A_Z$ be the linear map whose columns are these coefficients.
If every $z_i=c$, its columns are exactly those of $J_c$.
Uniform $C^4$ bounds and $|z_i-c|\le\rho$ therefore give
\[
 \|A_Z-J_c\|\le C\rho.
\]
Indeed, the degree-$k$ coefficient in the $k$th summand differs
from $D^kw(c)/k!$ by $O(\rho)$, while each lower-degree coefficient
in that summand contains the appropriate positive power of $\rho$.
Choosing $\rho$ small gives a uniform positive lower bound on
$\sigma(A_Z)$. Interpolation now gives the exact identity
\[
 \sum_i a_iw(z_i)=A_Z\left(\sum_i a_i(z_i-c)^\alpha\right)_{|\alpha|\le3},
\]
and hence \eqref{eq:fourmoment}.

For the compressed map $w_m=P_mw$, the $C^4$ bounds do not increase.
Moreover $P_mJ_c\to J_c$ uniformly in $c$, since the finitely many
jet columns range over compact subsets of $\Y$. Thus their least
singular values have one positive lower bound for all large $m$.
The same coefficient argument proves the final assertion.
\end{proof}

\begin{lemma}[Stability of positive two-atom moments]\label{lem:twoatom}
Let $\mu=a_1\delta_{x_1}+a_2\delta_{x_2}$ and
$\widetilde\mu=\widetilde a_1\delta_{y_1}+
\widetilde a_2\delta_{y_2}$, where all nodes lie in a fixed bounded
set, $a_j,\widetilde a_j\in[a_-,a_+]$ with $a_->0$, and each
within-pair separation is at least $0<h\le1$. Set
\[
 M_k(\mu)=\sum_j a_jx_j^{\otimes k},\qquad
 \eta=\max_{0\le k\le3}\|M_k(\mu)-M_k(\widetilde\mu)\|.
\]
There is one permutation $\pi\in\mathfrak S_2$ such that
\begin{equation}\label{eq:twoatomstability}
 \|X-Y_\pi\|\le Ch^{-2}\eta,\qquad
 \|a-\widetilde a_\pi\|\le Ch^{-3}\eta.
\end{equation}
The constant depends only on the bounded node set and the amplitude
interval.
\end{lemma}
\begin{proof}
From the raw moments form the total mass $A$, centroid $c$, covariance
$S$ and third central tensor $T$:
\[
 A=M_0,\quad c=M_1/A,\quad
 S=\sum_j a_j(x_j-c)^{\otimes2},\quad
 T=\sum_j a_j(x_j-c)^{\otimes3}.
\]
The differences of each of these quantities between the two measures
are at most $C\eta$: they are smooth functions of the first four
moments on a bounded set with $A\ge2a_->0$.
Put $t=a_1/A$, $q=1-t$, $r=|x_1-x_2|$ and
$e=(x_1-x_2)/r$. Then
\begin{equation}\label{eq:centralatom}
 \begin{gathered}
 x_1=c+qr e,\qquad x_2=c-tr e,\qquad
 t,q\in[\tau,1-\tau],\quad \tau>0,\\
 S=s\,ee^T,\quad s=Atqr^2,\qquad
 T=Atq(q-t)r^3e^{\otimes3}.
 \end{gathered}
\end{equation}
Use tildes for the second measure. If $\eta\ge\varepsilon h^2$
for a fixed small $\varepsilon>0$, boundedness of all parameters
proves \eqref{eq:twoatomstability}, with any permutation and a
sufficiently large constant. We may thus assume
$\eta<\varepsilon h^2$.

Since $s\asymp r^2$ and $r\ge h$, choosing $\varepsilon$ small makes
$s,\widetilde s$ comparable and hence $r,\widetilde r$ comparable.
Set $r_* =\min(r,\widetilde r)\ge h$. Normalize the rank-one
covariances and, by relabeling the second measure, choose
$e\cdot\widetilde e\ge0$. The rank-one projection identity gives
\begin{equation}\label{eq:atomdirection}
 |e-\widetilde e|\le C\|ee^T-\widetilde e\widetilde e^T\|
                    \le C\eta/r_*^2.
\end{equation}
The normalized third moment is
\begin{equation}\label{eq:atomskew}
 z=\frac{\sqrt A\,T[e,e,e]}{s^{3/2}}
   =\frac{1-2t}{\sqrt{t(1-t)}}=:f(t),\qquad
 f'(t)=-\frac1{2[t(1-t)]^{3/2}}.
\end{equation}
Using $T=O(r^3)$, the moment differences, comparable radii and
\eqref{eq:atomdirection} gives
\[
 |z-\widetilde z|\le C\eta/r_*^3,\qquad
 |t-\widetilde t|\le C\eta/r_*^3.
\]
For example, the direction error in the normalized tensor is
$O(\eta/r_*^2)$ and the denominator error is also
$O(\eta/r_*^2)$; both are bounded by $C\eta/r_*^3$ because the
radii have a fixed upper bound. The numerator error is directly
$O(\eta/r_*^3)$. The second assertion follows from $|f'|\ge4$.

Finally $r=(s/[At(1-t)])^{1/2}$, whence
\[
 |r-\widetilde r|\le
 C\bigl(\eta/r_*+r_*\eta+r_*|t-\widetilde t|\bigr)
 \le C\eta/r_*^2.
\]
Substitution into \eqref{eq:centralatom}, using
$|c-\widetilde c|\le C\eta$, proves the position bound.
The identities $a_1=At$, $a_2=A(1-t)$ prove the amplitude bound.
Both use the same orientation and hence the same permutation.
\end{proof}

\subsection{A global estimate with fixed finite observations}
\begin{theorem}[Uniform two-site inverse rates]\label{thm:sharptwoupper}
For any fixed complete real boundary basis there are $m_0$, $C>0$
and $0<h_0\le1$, independent of separation, such that for
$0<h\le h_0$ and $p,q\in\mathcal P_2(h)$ there is a permutation
$\pi$ with
\begin{equation}\label{eq:twosharpupper}
 \|X-Y_\pi\|\le Ch^{-2}D_{m_0}(p,q),\qquad
 \|b-\widetilde b_\pi\|\le Ch^{-3}D_{m_0}(p,q),
\end{equation}
where
$D_m(p,q)=\|P_m(F(p)-F(q))P_m\|$.
In particular $d_2(p,q)\le C'h^{-3}\|F(p)-F(q)\|$.
The same conclusions hold in the strengths $\alpha_j=1/b_j$.
The inputs and tests defining $P_{m_0}FP_{m_0}$ are fixed before
the positions, strengths and separation are known.
\end{theorem}
\begin{proof}
We first give a local estimate near an arbitrary collision centre $c$.
Let $a=w(c)/\|w(c)\|^2$, so that
$H(x)=\langle w(x),a\rangle$ satisfies $H(c)=1$.
For $p$ near that collision put
\[
 q_p=Q_p^{-1}(H(x_1),H(x_2))^T.
\]
At $x_1=x_2=c$, the Sherman--Morrison formula gives
\begin{equation}\label{eq:positivecharges}
 (q_p)_j=\frac{\alpha_j}{1+g_\lambda(\alpha_1+\alpha_2)}>0.
\end{equation}
Compactness of the $b$ interval and continuity therefore place both
charges in one interval $[a_-,a_+]\Subset(0,\infty)$ whenever the
sites are sufficiently close to $c$. These bounds are uniform in
the physical strengths.

For two such configurations the response column difference is
\[
 (F(p)-F(q))a
 =-\sum_j(q_p)_jw(x_j)+\sum_j(q_q)_jw(y_j).
\]
Lemma~\ref{lem:fourmoment} bounds its first four signed moments by
$C\|(F(p)-F(q))a\|$. Translation between monomials based at $c$
and raw tensor moments is uniformly bounded. Lemma~\ref{lem:twoatom}
then gives the position bound with $h^{-2}$ and the effective-charge
bound with $h^{-3}$. The physical inverse strengths are recovered
by the exact formula
\begin{equation}\label{eq:recoverbcharges}
 b_j=\frac{H(x_j)-\sum_{\ell=1}^2
                 G_\lambda(x_j-x_\ell)(q_p)_\ell}{(q_p)_j}.
\end{equation}
This is uniformly Lipschitz in the nodes and positive charges on
the sets under consideration. Indeed $H$ has bounded first derivatives,
and $G_\lambda$ is $C^1$ with bounded derivative, since
\[
 \int_{\R^2}\frac{|\xi|}{|\xi|^4+\lambda}\dd\xi<\infty.
\]
The denominator in \eqref{eq:recoverbcharges} stays uniformly
positive. Thus the same $h^{-3}$ bound holds for $b$.

This argument also works, with constants independent of large $m$,
for compressed data. Use $a_m=P_ma$, incoming field
$H_m(x)=\langle w(x),a_m\rangle$ and columns $w_m=P_mw$.
The incoming fields converge to $H$ in $C^1$, the charges retain
the bounds in \eqref{eq:positivecharges}, and
Lemma~\ref{lem:fourmoment} is uniform for $w_m$. The measured
column is exactly $P_m(F(p)-F(q))P_ma$, of norm at most
$\|a\|D_m(p,q)$. Formula~\eqref{eq:recoverbcharges} uses $H_m$
and retains its uniform Lipschitz bound. We have proved
\eqref{eq:twosharpupper} locally near every one-site collision,
uniformly for all sufficiently large $m$.

We now prove that one $m_0$ works globally. Let the labeled compact
parameter closure allow $x_1=x_2$. Its response remains continuous
because $Q\ge b_-I$. At a collision the two scalar forms merge
to the single strength $\alpha_1+\alpha_2$. The image of this
closure is compact in operator norm and consists of compact
operators. Consequently $P_mFP_m\to F$ uniformly on the closure.

Define the weighted relabeling discrepancy
\[
 \rho_h(p,q)=\min_\pi
 \bigl(h^2\|X-Y_\pi\|+h^3\|b-\widetilde b_\pi\|\bigr).
\]
It suffices to prove $\rho_h\le CD_{m_0}$. If no such fixed
$m_0,C$ existed for $0<h\le h_0$, we could choose
$p_m,q_m\in\mathcal P_2(h_m)$ with
\[
 \rho_{h_m}(p_m,q_m)>mD_m(p_m,q_m).
\]
The left side is bounded, so $D_m(p_m,q_m)\to0$. Uniform
compression convergence implies
$\|F(p_m)-F(q_m)\|\to0$. Take a subsequence with limiting
labeled configurations $p_0,q_0$. Their merged responses agree.
Proposition~\ref{thm:sites}, applied after merging coincident positive
sites, shows that their merged configurations agree.

If that common configuration has two distinct sites, the strengths
also agree after relabeling. Thus $p_m,q_m$ tend to the same
regular parameter. The derivative is injective by
Proposition~\ref{thm:sites}. Its least singular value at that
parameter is positive. Compression of that derivative converges in
operator norm on its finite-dimensional domain. A uniform $C^1$
expansion in a small separated coordinate ball therefore gives
\[
 D_m(p_m,q_m)\ge c\|p_m-q_m\|
\]
for all large $m$, contradicting the displayed inequality because
$h_m\le1$.

If the merged configuration has one site $c$, all four nodes tend
to $c$. The compressed local estimate already proved applies for
large $m$ and gives $\rho_{h_m}\le CD_m$, again a contradiction.
These alternatives exhaust the compact closure and prove the
uniform finite-compression estimate. Since $h\le1$, it implies
the joint $h^{-3}$ bound. The reciprocal map is bi-Lipschitz on
$[b_-,b_+]$, proving the assertion for $\alpha$.
\end{proof}

The proof uses an incoming field associated with a collision centre
only to establish an estimate. The actual experiment records the
fixed finite matrix $P_{m_0}FP_{m_0}$; its linear combinations
supply every column used in the argument. No knowledge of the centre
is needed to choose the inputs. This theorem asserts existence of a
finite matrix size, not an optimal count or a quantitative bound on $m_0$.

\subsection{Matching physical lower bounds}
\begin{theorem}[Sharpness in the nonlinear physical model]
\label{thm:sharptwolower}
Suppose $K$ contains a disk. For all sufficiently small $h$, there
are $p_h,q_h\in\mathcal P_2(h)$ such that
\begin{equation}\label{eq:twolowerfamily}
 \|F(p_h)-F(q_h)\|\asymp h^3,\qquad
 \min_\pi\|X_h-Y_{h,\pi}\|\asymp h,\qquad
 \min_\pi\|b_h-\widetilde b_{h,\pi}\|\asymp1.
\end{equation}
Consequently the best joint inverse Lipschitz constant on
$\mathcal P_2(h)$ satisfies $C_{\rm joint}(h)\asymp h^{-3}$.
The best position constant is $C_{\rm pos}(h)\asymp h^{-2}$,
and the best inverse-strength constant is
$C_b(h)\asymp h^{-3}$. The same powers hold for the fixed finite
compression in Theorem~\ref{thm:sharptwoupper} and for the strengths
$\alpha=1/b$.
\end{theorem}
\begin{proof}
Choose a centre $c$ in the disk, a unit vector $e$, and
$b_0\in(b_-,b_+)$. Set
\[
 g=g_\lambda,\qquad L=b_0/2,\qquad S=(g+L)^{-1}.
\]
For $t$ in a fixed small interval about $1/2$ put $q=1-t$ and
\begin{equation}\label{eq:physicalmatchingfamily}
 \begin{gathered}
 r_t=\frac{h}{2\sqrt{tq}},\qquad
 x_1=c+qr_te,\qquad x_2=c-tr_te,\qquad
 V=tqr_t^2=h^2/4,\\
 c_t=G_\lambda(r_t),\quad \delta_t=g-c_t,\qquad
 b_1=\frac{S^{-1}-c_t}{t}-\delta_t,\qquad
 b_2=\frac{S^{-1}-c_t}{q}-\delta_t.
 \end{gathered}
\end{equation}
Here $G_\lambda(r)$ denotes its radial value. The separation
$r_t\ge h$, and $b_1\to L/t$, $b_2\to L/q$ as $h\to0$.
Restricting the $t$ interval makes these limits and, for small $h$,
the actual strengths lie strictly in $[b_-,b_+]$. All nodes remain
in the chosen disk.

The following identities verify that this is a family of the
physical diagonal-coupling model. With $v=(t,q)^T$ and
$d=(1,-1)^T$, direct inversion gives
\begin{equation}\label{eq:physicalBdecomp}
 \begin{aligned}
 Q&=(L+\delta_t)\diag(1/t,1/q)+c_t\mathbf1\mathbf1^T,\\
 Q^{-1}&=Svv^T+\kappa_tdd^T,\qquad
 \kappa_t=\frac{tq}{L+\delta_t}.
 \end{aligned}
\end{equation}
In particular $\mathbf1^TQ^{-1}\mathbf1=S$ is exactly independent
of $t$. If $w_j=w(x_j)$ and $\overline w_t=tw_1+qw_2$, then
\begin{equation}\label{eq:physicalFdecomp}
 F_t=-S\overline w_t\overline w_t^*
             -\kappa_t(w_1-w_2)(w_1-w_2)^*.
\end{equation}

We also need to control the internal Green correction. Its Fourier
representation and $1-\cos u\le C\min(u^2,1)$ give, for small $r$,
\begin{equation}\label{eq:greenquadraticlog}
 0\le g-G_\lambda(r)
 \le C\int_0^\infty\frac{\min(r^2s^2,1)s}{s^4+\lambda}\dd s
 \le Cr^2\log(1/r).
\end{equation}
The last step splits the integral at a fixed scale and at $1/r$.
Thus $\delta_t=O(h^2\log(1/h))$, uniformly for the chosen $t$.

Use source-direction derivatives $w_e,w_{ee},w_{eee}$ at $c$,
and write $w=w(c)$. Taylor expansion in $\Y$ gives
\[
 \begin{aligned}
 \overline w_t&=w+\frac V2w_{ee}+\frac{m_3}6w_{eee}+O(h^4),\\
 w_1-w_2&=r_tw_e+\frac{q-t}{2}r_t^2w_{ee}+O(h^3),\\
 m_3&=tq(q-t)r_t^3=h^3\gamma(t),\qquad
 \gamma(t)=\frac{1-2t}{8\sqrt{tq}}.
 \end{aligned}
\]
Substitution into \eqref{eq:physicalFdecomp}, using
\eqref{eq:greenquadraticlog}, yields the operator-norm expansion
\begin{equation}\label{eq:physicalthirdexpansion}
 \begin{aligned}
 F_t={}&-Sww^*-\frac{SV}{2}(ww_{ee}^*+w_{ee}w^*)
              -\frac VLw_ew_e^*-m_3\mathcal B
              +O(h^4\log(1/h)),\\
 \mathcal B={}&\frac S6(ww_{eee}^*+w_{eee}w^*)
              +\frac1{2L}(w_ew_{ee}^*+w_{ee}w_e^*).
 \end{aligned}
\end{equation}
All terms through degree two are independent of $t$.
The coefficient $\mathcal B$ is nonzero. Indeed let $P$ project
orthogonally onto the complement of
$\operatorname{span}\{w,w_e,w_{ee}\}$. Then
\[
 P\mathcal B w=\frac{S\|w\|^2}{6}Pw_{eee}\ne0
\]
by the order-three jet independence in Lemma~\ref{lem:arc}.

Compare $t_0=1/2$ with a fixed $t_1\ne1/2$ in the chosen interval.
Their response difference is
\[
 F_{t_1}-F_{t_0}=-h^3\gamma(t_1)\mathcal B
                      +O(h^4\log(1/h)),
\]
so its norm is asymptotic to
$|\gamma(t_1)|\,\|\mathcal B\|h^3>0$.
At $t_0$, both inverse strengths tend to $b_0$; at $t_1$, they
tend to $L/t_1,L/(1-t_1)$, a different unordered pair. Their
strength distance therefore tends to a positive constant.
The optimal positional matching preserves the direction of the
ordered pair, and
\[
 \frac1h\min_\pi\|X_{t_1}-X_{t_0,\pi}\|
 =\frac12\left[
    \left(\sqrt{\frac{1-t_1}{t_1}}-1\right)^2+
    \left(\sqrt{\frac{t_1}{1-t_1}}-1\right)^2
             \right]^{1/2}>0.
\]
For example, the direct-minus-exchanged squared matching cost is
$-2r_{t_1}r_{t_0}<0$. This proves
\eqref{eq:twolowerfamily}. It gives the three necessary lower
powers for the full response, while
Theorem~\ref{thm:sharptwoupper} supplies their upper bounds.
Compression cannot increase a response norm, so the same lower
bounds hold for its fixed finite data. Its upper bounds are already
proved. Finally the reciprocal map preserves these powers.
\end{proof}

\begin{corollary}[Sharp deterministic noise threshold]\label{cor:noiseboundary}
Suppose $K$ contains a disk and let $\eta(h)\ge0$ bound additive operator-norm noise. Uniform recovery on $\mathcal P_2(h)$ with position error $o(h)$ and inverse-strength error $o(1)$ under a common relabeling is possible if and only if $\eta(h)=o(h^3)$. This holds for the full response and for the fixed compression of Theorem~\ref{thm:sharptwoupper}. The same criterion holds for physical strengths.
\end{corollary}
\begin{proof}
A residual minimizer over the compact class has data within $2\eta(h)$ of the true response. Theorem~\ref{thm:sharptwoupper} gives errors at most $2Ch^{-2}\eta(h)$ and $2Ch^{-3}\eta(h)$ under one matching. A full observation can be compressed with operator norm one.

Conversely, if $\eta(h)/h^3$ fails to tend to zero, then $\eta(h)\ge\epsilon h^3$ on a sequence $h\downarrow0$ for some $\epsilon>0$. In \eqref{eq:physicalmatchingfamily} take $t_0=1/2$ and fix $t_1\ne t_0$ sufficiently close to $t_0$ that $|\gamma(t_1)|\,\|\mathcal B\|<2\epsilon$. By \eqref{eq:physicalthirdexpansion}, the midpoint of the two exact full responses lies within $\eta(h)$ of both along the sequence, for all sufficiently small $h$. Its compression has the same property. The two unordered position pairs have distance at least $c(t_1)h$ and the two unordered strength pairs have distance at least $c'(t_1)>0$ by Theorem~\ref{thm:sharptwolower} and its construction. In each quotient metric the triangle inequality forces any estimate at the common datum to incur at least half that distance for one of the two configurations. Reciprocation is bi-Lipschitz on the prescribed positive interval.
\end{proof}

\section*{Declarations}

\subsection*{Funding}
This research was funded by the Committee of Science of the Ministry of
Science and Higher Education of the Republic of Kazakhstan (Grant
Nos.~BR31714735 and AP32319153).

\subsection*{Conflict of interest}
The author declares that there is no conflict of interest.

\subsection*{Data availability}
No datasets were generated or analysed during the present study.

\subsection*{Code availability}
Not applicable.

\subsection*{Ethics approval}
Not applicable.

\subsection*{Declaration of generative AI use}
The author used ChatGPT (OpenAI) and Claude (Anthropic) for language
editing. ChatGPT assisted in revising two arguments: the use of
invariance of domain for the exchange loop in
Proposition~\ref{prop:topologicalcount} and the noise lower bound in
Corollary~\ref{cor:noiseboundary}. Claude checked the proofs, carried out
symbolic and numerical checks of selected identities and estimates,
verified the bibliographic data and drafted Remark~\ref{rem:shiftsize}.
The author reviewed and edited all such output and takes full
responsibility for the content of the publication.

\end{document}